\documentclass[12pt]{article}
\usepackage{amscd}
\usepackage{latexsym}
\usepackage{mathrsfs}
\usepackage{pst-all,multido,ifthen}
\usepackage{amsmath}
\usepackage{amssymb}
\usepackage[all]{xy}
\usepackage[english]{babel}
\newtheorem{thm}{Theorem}
\newtheorem{df}{Definition}
\newtheorem{prop}{Proposition}
\newtheorem{cor}{Corollary}
\newtheorem{lemma}{Lemma}

\newtheorem{fact}{Fact}

\newtheorem{remark}{Remark}
\def\endproof{$\hfill \square$}
\def\Spec{\text{Spec\;}}
\def\Hom{\text{Hom\;}}
\def\perf{\text{perf\;}}

\begin{document}
\title{Purity of Crystalline Strata}
\author{Jinghao Li and Adrian Vasiu}
\maketitle

\noindent
{\bf ABSTRACT.} Let $p$ be a prime. Let $n\in\mathbb N^{\ast}$. Let $\mathcal C$ be an $F^n$-crystal over a locally noetherian $\mathbb F_p$-scheme $S$. Let $(a,b)\in\mathbb N^2$. We show that the reduced locally closed subscheme of $S$ whose points are exactly those $x\in S$ such that $(a,b)$ is a break point of the Newton polygon of the fiber $\mathcal C_x$ of $\mathcal C$ at $x$ is pure in $S$, i.e., it is an affine $S$-scheme. This result refines and reobtains previous results of de Jong--Oort, Vasiu, and Yang. As an application, we show that for all $m\in \mathbb N$ the reduced locally closed subscheme of $S$ whose points are exactly those $x\in S$ for which the $p$-rank of $\mathcal C_x$ is $m$ is pure in $S$; the case $n=1$ was previously obtained by Deligne (unpublished) and the general case $n\ge 1$ refines and reobtains a result of Zink.

\bigskip\noindent
{\bf KEY WORDS:} $\mathbb F_p$-scheme, $F$-crystal, Newton polygon, $p$-rank, purity.

\bigskip\noindent
{\bf MSC 2010:} 11G10, 11G18, 14F30, 14G35, 14K10, 14K99, 14L05, and 14L15

\section{Introduction}\label{S1}
For a reduced locally closed subscheme $Z$ of a locally noetherian scheme $Y$, let $\bar Z$ be the schematic closure of $Z$ in $Y$. We recall from \cite{NVW}, Definition 1.1 that $Z$ is called {\it pure} in $Y$ if it is an affine $Y$-scheme. The paper \cite{NVW} also uses a weaker variant of this purity which in \cite{L} is called {\it weakly pure}: we say $Z$ is weakly pure in $Y$ if each non-empty irreducible component of the complement $\bar Z-Z$ is of pure codimension $1$ in $\bar Z$. It is well-known that if $Z$ is pure in $Y$, then $Z$ is also weakly pure in $Y$ (for instance, cf. Proposition \ref{P3} of Subsection \ref{S44}).

Let $n$ and $r$ be natural numbers. Let $p$ be a prime. Let $S$ be a locally noetherian $\mathbb F_p$-scheme. Let $\Phi_S:S\to S$ be the Frobenius endomorphism of $S$. Let $\mathcal M$ be a {\it crystal} of the gross absolute crystalline site $CRIS(S/\Spec (\mathbb Z_p))$ introduced in \cite{B}, Chapter III, Example 1.1.3 and Definition 4.1.1 in locally free $\underline{\mathcal O}_{S/\Spec (\mathbb Z_p)}$-modules of rank $r$. We assume that we have an {\it isogeny} $\phi_{\mathcal M}:(\Phi_S^n)^*(\mathcal M)\to\mathcal M$; thus  the pair $\mathcal C=(\mathcal M, \phi_{\mathcal M})$ is an $F^n$-crystal of $CRIS(S/\Spec (\mathbb Z_p))$. If the $\mathbb F_p$-scheme $S=\Spec A$ is affine, then the pair $\mathcal C=(\mathcal M, \phi_{\mathcal M})$ is canonically identified with a $\sigma^n-F$-crystal on $A$ in the sense of \cite{K}, Subsection (2.1). 

Let $\nu:[0,r]\to [0,\infty)$ be a {\it Newton polygon}, i.e., a nondecreasing piecewise linear continuous function such that $\nu(0)=0$ and the coordinates of all its {\it break points} are natural numbers. For $x\in S$, let $\nu_x$ be the {\it Newton polygon} of the fiber $\mathcal C_x$ of $\mathcal C$ at $x$. Let $S_{\nu}$ be the reduced locally closed subscheme of $S$ whose points are exactly those $x\in S$ such that we have $\nu_x=\nu$, cf. Grothendieck--Katz Theorem (see \cite{K}, Corollary 2.3.2); if non-empty, $S_{\nu}$ is a {\it stratum} of the Newton polygon {\it stratification} of $S$ defined by $\mathcal C$.

Let $a,b\in\mathbb N$ be such that $0\le a\le r$. Let $T=T_{(a,b)}(\mathcal C)$ be the reduced locally closed subscheme of $S$ whose points are those $x\in S$ such that $(a,b)$ is a break point of $\nu_x$. The end break point $(r,\nu_x(r))$ remains constant under specializations of $x\in S$. Thus locally in the Zariski topology of $S$, we can assume that there exists $d\in\mathbb N$ such that for all $x\in S$ we have $\nu_x(r)=d$ and this implies that $T$ is the reduced locally closed subscheme of $S$ which is a finite union $\bigcup_{\nu\in N_{r,d,a,b}} S_{\nu}$ of Newton polygon strata $S_{\nu}$ indexed by the set $N_{r,d,a,b}$ of all Newton polygons $\nu:[0,r]\to  [0,\infty)$ with the two properties that $\nu(r)=d$ and $(a,b)$ is a break point of $\nu$. 

It is known that $T$ is weakly pure in $S$, cf. \cite{Y}, Theorem 1.1. It is also known that $S_{\nu}$ is pure in $S$, cf. \cite{V1}, Main Theorem B. This last result implies the celebrated result of de Jong--Oort which asserts that $S_{\nu}$ is weakly pure in $S$, cf. \cite{dJO}, Theorem 4.1. Strictly speaking, the references of this paragraph work with $n=1$ but their proofs apply to all $n\in\mathbb N^{\ast}$. 

In general, a finite union of locally closed subschemes of $S$ which are pure in $S$ is not pure in $S$. Therefore the following purity result which refines and reobtains the mentioned results of de Jong--Oort, Vasiu, and Yang comes as a surprise. 

\begin{thm}\label{T1}
With the above notations, $T$ is pure in $S$.
\end{thm}

In Section \ref{S2} we gather few preliminary steps that are required to prove Theorem \ref{T1} in Section \ref{S3}. We have the following two direct consequences of Theorem \ref{T1}, the first one for $n=1$ just reobtains \cite{V1}, Main Theorem B in the locally noetherian case. 

\begin{cor}\label{C1}
Each Newton polygon stratum $S_{\nu}$ is pure in $S$.
\end{cor}

The $p$-rank $\chi(x)$ of $\mathcal C_x$ is the multiplicity of the Newton polygon slope $0$ of $\nu_x$. Equivalently, $\chi(x)$ is the unique natural number such that $(0,0)$ and $(\chi(x),0)$ are the only break points of $\nu_x$ on the horizontal axis (i.e., which have the second coordinate $0$). 

\begin{cor}\label{C2}
Let $m\in \mathbb N$. We consider the reduced locally closed subscheme $S_m$ of $S$ whose points are exactly those $x\in S$ such that the $p$-rank $\chi(x)$ of $\mathcal C_x$ is $m$.  Then $S_m$ is pure in $S$.
\end{cor}

If $m>0$, then we have $S_m=T_{(m,0)}(\mathcal C)$ and if  $m=0$, then we have $S_0=T_{(1,0)}(\mathcal C\oplus \mathcal E_0)$ where $\mathcal E_0$ is the pull back to $S$ of the $F^n$-crystal over $\Spec (\mathbb F_p)$ of rank $1$ and Newton polygon slope $0$ which has a Frobenius invariant global section; therefore, regardless of what $m$ is, Corollary \ref{C2} follows from Theorem \ref{T1}.

For $n=1$ Corollary \ref{C2} was first obtained by Deligne and more recently by Vasiu and Li (see \cite{D}, \cite{V3}, and \cite{L}). Corollary \ref{C2} also refines and reobtains a prior result of Zink which asserts that $S_m$ is weakly pure in $S$ (see \cite{Z}, Proposition 5). 

In Section \ref{S4} we first follow \cite{L} to show that Corollary \ref{C1} follows directly from Theorem \ref{T1} and then we follow \cite{V3} to include a second proof of Corollary \ref{C2} in the more general context provided by a functorial version of the {\it Artin--Schreier stratifications} introduced in \cite{V2}, Definition 2.4.2 which is simpler, does not rely on Theorem \ref{T1}, and is based on Theorem \ref{T2} of Subsection \ref{S42}. 

Theorem \ref{T1} is due to the first author, cf. \cite{L}. While the proof of \cite{Y}, Theorem 1.1 follows the proof of \cite{dJO}, Theorem 4.1, the proof of Theorem \ref{T1} presented follows \cite{L} and thus the proofs of \cite{V1}, Main Theorem B and Theorem 6.1. It is known (cf. \cite{NVW}, Example 7.1) that in general $S_m$ is not strongly pure in $S$ in the sense of \cite{NVW}, Definition 7.1 and therefore Theorem \ref{T1} and Corollary \ref{C2} cannot be improved in general (i.e., are optimal).

We refer to $T_{(a,b)}(\mathcal C)$, $S_{\nu}$, and $S_m$ as crystalline strata of $S$ associated to $\mathcal C$ and certain (basic) discrete invariants of $F^n$-crystals. Cases of non-discrete invariants stemming from isomorphism classes are also studied in the literature (for instance, see \cite{V1}, Subsection 5.3 and \cite{NVW}, Theorem 1.2 and Corollary 1.5). Crystalline strata have applications to the study in positive characteristic of different moduli spaces and schemes such as special fibers of Shimura varieties of Hodge type (for instance, see \cite{V1} and \cite{NVW}). 

\section{Standard reduction steps}\label{S2}

The above notations $p$, $S$, $\Phi_S$, $\bar Z$, $n$, $r$, $\mathcal C=(\mathcal M,\phi_{\mathcal M})$, $\mathcal C_x$, $\nu_x$, $(a,b)\in \mathbb N^2$, $T=T_{(a,b)}(\mathcal C)$, $S_{\nu}$, $m$, $S_m$, $\chi(x)$, and $\mathcal E_0$ will be used throughout the paper. For a fixed Newton polygon $\nu$, let $S_{\ge \nu}$ be the reduced closed subscheme of $S$ whose points are exactly those $x\in S$ such that the Newton polygon $\nu_x$ is above $\nu$, cf. \cite{K}, Corollary 2.3.2. 

In what follows by an \'etale cover we mean a surjective finite \'etale morphism of schemes. For basic properties of excellent rings we refer to \cite{M}, Chapter 13. If $V\to Y$ is a morphism of $\mathbb F_p$-schemes and if $\mathcal F$ (or $\mathcal F_Y$) is an $F^n$-crystal over $Y$, let $\mathcal F_V$ be the pull back of $\mathcal F$ (or $\mathcal F_Y$) to an $F^n$-crystal over $V$, i.e., of $CRIS(V/\Spec (\mathbb Z_p))$. Let $k(y)$ be the residue field of a point $y\in Y$. If $V=\Spec (k(y))\to Y$ is the natural morphism, then we denote $\mathcal F_V=\mathcal F_{\Spec (k(y))}$ simply by $\mathcal F_y$ (the fiber of $\mathcal F$ at $y$).

For an $\mathbb F_p$-algebra $R$, let $W(R)$ be the ring of $p$-typical Witt vectors with coefficients in $R$. Let $\mathbb W(R)=(\Spec R,\Spec (W(R)),\text{can})$ be the thickening in which `$\text{can}$' stands for the canonical divided power structure of the kernel of the epimorphism $W(R)\to W_1(R)=R$. For $s\in\mathbb N^{\ast}$, let $W_s(R)$ be the ring of $p$-typical Witt vectors of length $s$ with coefficients in $R$. Let $\mathbb W_s(R)=(\Spec R,\Spec (W_s(R)),\text{can})$ be the thickening defined naturally by $\mathbb W(R)$. Let $\Phi_R$ be the Frobenius endomorphism of either $W(R)$ or $W_s(R)$.

The property of a reduced locally closed subscheme being pure in $S$ is local for the faithfully flat topology of $S$, and thus until the end we will also assume that $S=\Spec A$ is an affine $\mathbb F_p$-scheme and that there exists $d\in\mathbb N$ such that for all $x\in S$ we have $\nu_x(r)=d$. As the scheme $S$ is locally noetherian and affine, it is noetherian. To prove Theorem \ref{T1}, we have to prove that $T$ is an affine scheme.

\subsection{Some abelian categories}\label{S21}

Let $\mathcal M(W_s(R))$ be the abelian category whose objects are pairs $(O,\phi_O)$ comprising from a $W_s(R)$-module $O$ and a $\Phi_R^n$-linear endomorphism $\phi_O:O\to O$ (i.e., $\phi_O$ is additive and for all $z\in O$ and $\sigma\in W_s(R)$ we have $\phi_O(\sigma z)=\Phi_R^n(\sigma)\phi_O(z)$) and whose morphisms $f:(O_1,\phi_{O_1})\to (O_2,\phi_{O_2})$ are $W_s(R)$-linear maps $f:O_1\to O_2$ satisfying $f\circ\phi_{O_1}=\phi_{O_2}\circ f$. If $t\in\{0,\ldots,s-1\}$, then by a {\it quasi-isogeny} of $\mathcal M(W_s(R))$ whose cokernel is annihilated by $p^t$ we mean a morphism $f:(O_1,\phi_{O_1})\to (O_2,\phi_{O_2})$ of $\mathcal M(W_s(R))$ which has the following two properties: (i) both $O_1$ and $O_2$ are projective $W_s(R)$-modules which locally in the Zariski topology of $\Spec (W_s(R))$ have the same positive rank, and (ii) the cokernel $O_2/f(O_1)$ is annihilated by $p^t$. An object $(O,\phi_O)$  of $\mathcal M(W_s(R))$ is called {\it divisible} by $t\in\{1,\ldots,s-1\}$ if $O$ is a projective $W_s(R)$-module such that $\text{Im}(\phi_O)\subseteq p^tO_2$. 

For $l\in\mathbb N^{\ast}$ we have a natural functor $\mathcal M(W_{s+l}(R))\to \mathcal M(W_s(R))$ to be referred by abuse of language as the reduction modulo $p^s$ functor.

If $Y$ is a $\Spec (\mathbb F_p)$-scheme, in a similar way we define the scheme $W_s(Y)$, its Frobenius endomorphism $\Phi_Y$, and the abelian category $\mathcal M(W_s(Y))$ and speak about quasi-isogenies of $\mathcal M(W_s(Y))$ whose cokernels are annihilated by $p^t$ with $t\in\{0,\ldots,s-1\}$, about objects of $\mathcal M(W_s(Y))$ divisible by $t\in\{1,\ldots,s-1\}$, and about reduction modulo $p^s$ functors $\mathcal M(W_{s+l}(Y))\to \mathcal M(W_s(Y))$. We have canonical identifications $\mathcal M(W_s(R))=\mathcal M(W_s(\Spec R))$.

For homomorphisms $R\to R_1$ and morphisms $Y_1\to Y$ we have natural pull back functors $\mathcal M(W_s(R))\to \mathcal M(W_s(R_1))$ and $\mathcal M(W_s(Y))\to \mathcal M(W_s(Y_1))$. 

To prove that $T$ is an affine scheme, we can also assume that the {\it evaluation} $M$ of $\mathcal M$ at the thickening $\mathbb W_1(A)$ is a free $A$-module of rank $r$. The evaluation of $\phi_{\mathcal M}$ at this thickening is a $\Phi_A^n$-linear endomorphism $\phi_M:M\to M$. 

In what follows we will apply twice the following elementary general fact which can be also deduced easily from the elementary divisor theorem.

\begin{fact}\label{F1}
Let $D$ be a discrete valuation ring and let $\pi\in D$ be a uniformizer of it. Let $s,t\in\mathbb N$ be such that $s>t$. Let $D_s=D/(\pi^s)$. Let $g_s:D_s^r\to D_s^r$ be a $D_s$-linear endomorphism such that its cokernel is annihilated by $\pi^t$. Then for each $x\in D_s^r-\pi D_s^r$, we have $g_s(x)\in D_s^r-\pi^{t+1}D_s$.
\end{fact}

\noindent
{\bf Proof:} Let $g:D^r\to D^r$ be a $D$-linear endomorphism which lifts $g_s$. Let $E=\text{Im}(g)+\pi^sD^r$ (one can easily check that $E=\text{Im}(g)$ but we will not stop to argue this). It is a free $D$-module of rank $r$ which (as $\pi^t\text{Coker}(g_s)=0$) contains $\pi^tD^r$.  Thus $\pi^sD^r\subseteq pE$ and therefore $\text{Im}(g)$ surjects onto the $D_1$-vector space $E/\pi E$ of rank $r$. Hence a $D_s$-basis of $D_s^r$ maps via $g$ to a $D_1$-basis of $E$. From this and the fact that $\pi^{t+1}D^r\subseteq \pi E$ we get that no element of a $D_s$-basis of $D_s^r$ is mapped by $g$ to $\pi^{t+1}D^r$. Thus the fact holds. 

\subsection{On $(a,b)$}\label{S22}

If $(a,b)$ is $(0,0)$ or $(r,d)$, then $T=S$. If $a=0$ and $b>0$ or if $a=r$ and $b\neq d$, then $T=\emptyset$. Thus, to prove that $T$ is an affine scheme we can assume that $1\le a\le r-1$.

\begin{lemma}\label{L1}
Let $k$ be a field of characteristic $p$. Let $\nu:[0,r]\to [0,\infty)$ be the Newton polygon of an $F^n$-crystal $\mathcal F$ over $k$ of rank $r$. Let $a,b\in\mathbb N$ be such that $1\le a\le r-1$. Then $(a,b)$ is a break point of $\nu$ if and only if $(1,b)$ is a break point of the Newton polygon $\bigwedge^a(\nu)$ of the $F^n$-crystal  over $k$ of rank $\binom{r}{a}$ which is the exterior power $\bigwedge^a(\mathcal F)$ of $\mathcal F$.
\end{lemma}

\noindent
{\bf Proof:} Let $\alpha_1\le\cdots\le\alpha_r$ be the Newton polygon slopes of $\nu$. Let $\beta_1\le\cdots\le\beta_{\binom{r}{a}}$ be the Newton polygon slopes of $\bigwedge^a(\nu)$. We have 
$$\beta_1=\sum_{i=1}^a \alpha_i\;\;\;\text{and}\;\;\; \beta_2=(\sum_{i=1}^{a-1} \alpha_i)+\alpha_{a+1}=\beta_1+\alpha_{a+1}-\alpha_a.$$
Thus $\beta_1<\beta_2$ if and only if $\alpha_a<\alpha_{a+1}$. Moreover, $(a,b)$ is a break point of $\nu$ if and only if we have $\alpha_a<\alpha_{a+1}$, and $(1,b)$ is a break point of the Newton polygon $\bigwedge^a(\nu)$ if and only if we have $\beta_1<\beta_2$. The lemma follows from the last two sentences.\endproof

\medskip
Based on Lemma \ref{L1}, to prove that $T$ is an affine scheme by replacing $\mathcal C$ with its exterior power $\bigwedge^a(\mathcal C)$ we can assume that $a=1$.  

\subsection{A description of $T$}\label{S23}

Let $q\in\mathbb N^{\ast}$ be such that for each $x\in S$ the Newton polygon slopes of the $F^{nq}$-crystal over $\Spec (k(x))$ which is the $q$-th iterate of $\mathcal C_x$ are all integers. For instance, as each Newton polygon slope of $\mathcal C_x$ is a rational number whose denominator is a natural number at most equal to $r$, we can take $q=r!$. Thus by replacing $n$ by $nq$ and $\mathcal C$ by its $q$-th iterate, we can assume that for each $x\in S$ the Newton polygon slopes of $\mathcal C_x$ are natural numbers. 

We consider the Newton polygon $\nu_1:[0,r]\to [0,\infty)$ whose graph is:

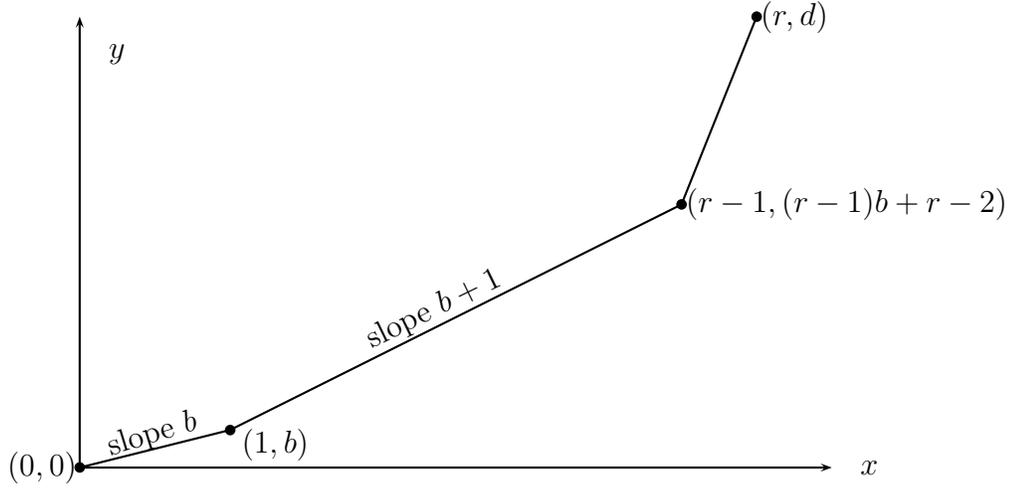
\begin{center}
\begin{pspicture}(10,7)
\psdots[dotsize=4pt](0,1)(2,1.5)(8,4.5)(9,7)
 \psline(0,1)(2,1.5)
\psline(2,1.5)(8,4.5) \psline(8,4.5)(9,7)
\psline{->}(0,1)(10,1)
\psline{->}(0,1)(0,7)
\rput(10.5,1){$x$}
\rput(0.5,6.5){$y$}
\rput[tr]{17}(1.5,1.8){\text{slope $b$}}
\rput[tr]{26}(5.5,3.7){\text{slope $b+1$}}
\rput(-0.5,1){$(0,0)$}
\rput(2.6,1.3){$(1,b)$}
\rput(10.2,4.5){$(r-1,(r-1)b+r-2)$}
\rput(9.5,7){$(r,d)$}
\end{pspicture}
\end{center}

If $x\in T$, then as all Newton polygon slopes of $\mathcal C_x$ are natural numbers, these Newton polygon slopes are $\alpha_1=b$, $\alpha_2\ge b+1$, $\alpha_{r-1}\ge b+1$, and $\alpha_r=d-\sum_{i=1}^{r-1} \alpha_i\ge b+1$. Therefore, if $x\in T$ then we have $x\in S_{\ge \nu_1}$. This implies that $T$ is a subscheme of the closed subscheme $S_{\ge \nu_1}$ of $S$. By replacing $S$ with $S_{\ge \nu_1}$ we can assume that $S=S_{\ge \nu_1}$. Thus $S$ is reduced.

If $r(b+1)>d$, then $S=S_{\ge \nu_1}=S_{\nu_1}=T$ and thus $T$ is affine. Thus we can assume that $r(b+1)\le d$ and therefore there exists a Newton polygon $\nu_2:[0,r]\to [0,\infty)$ whose graph is:

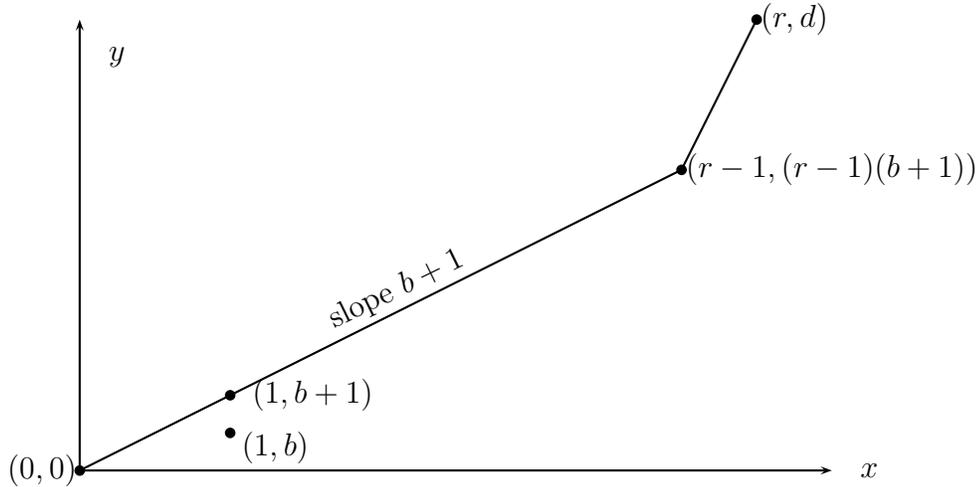
\begin{center}
\begin{pspicture}(10,7)
\psdots[dotsize=4pt](0,1)(2,1.5)(2,2)(8,5)(9,7)
\psline(0,1)(8,5)
\psline(8,5)(9,7)
\psline{->}(0,1)(10,1)
\psline{->}(0,1)(0,7)
\rput(10.5,1){$x$}
\rput(0.5,6.5){$y$}
\rput[tr]{26}(5,4){\text{slope $b+1$}}
\rput(-0.5,1){$(0,0)$}
\rput(3.1,2){$(1,b+1)$}
\rput(2.6,1.3){$(1,b)$}
\rput(10,5){$(r-1,(r-1)(b+1))$}
\rput(9.5,7){$(r,d)$}
\end{pspicture}
\end{center}

If $x\in S-T=S_{\ge \nu_1}-T$, then all Newton polygon slopes of $\mathcal C_x$ are natural numbers $\alpha_1\ge b+1$, $\alpha_2\ge b+1$, $\alpha_{r-1}\ge b+1$, and $\alpha_r=d-\sum_{i=1}^{r-1} \alpha_i\ge b+1$ and thus $\nu_x$ is above $\nu_2$. If $\nu_x$ is not above $\nu_2$, then as $\nu_x$ is above $\nu_1$ (as $S=S_{\ge\nu_1}$) we get that we have $\alpha_1=b$ and $\alpha_i\ge b+1$ for $i\in\{2,\ldots,r\}$. 

From the last two sentences we get that we have identities 
$$T=T_{(1,b)}=S - S_{\ge\nu_2}=S_{\ge\nu_1} - S_{\ge\nu_2}.$$
Thus, under all the above reduction steps, $T$ is an open subscheme of $S$.

\subsection{On $S$}\label{S24}

The statement that $T$ is an affine scheme is local in the faithfully flat topology of $S$ and therefore until the end of Section \ref{S3} we will assume that $A$ is a complete local reduced noetherian ring. Thus $A$ is also excellent and therefore its normalization in its ring of fractions is a finite product of normal complete local noetherian integral domains. Based on \cite{V1}, Lemma 2.9.2  which is a standard application of Chevalley's theorem of \cite{G1}, Chapter II, (6.7.1), to prove that $T$ is an affine scheme we can replace $A$ by one of the factors of the last product. Thus we can assume that $A$ is a normal complete local noetherian integral domain. We can also assume that $T$ is non-empty and therefore it is an open dense subscheme of $S$. Let $K$ be the field of factions of $A$ and let $\bar K$ be an algebraic closure of it. 

\section{Proof of Theorem \ref{T1}}\label{S3}

In this section we complete the proof of Theorem \ref{T1}, i.e., we prove that $T$ is an affine scheme when $a=1<r$, for each $x\in S$ all Newton polygon slopes of $\mathcal C_x$ are natural numbers, we have $S=S_{\ge\nu_1}=\Spec A$ with $A$ a normal complete local noetherian integral domain, and $T=T_{(1,b)}=S - S_{\ge\nu_2}$ is open dense in $S$. Let $\mathcal E_b=(\mathcal M_b,\phi_{\mathcal M_b})$ be the pull back to $S$ of the $F^n$-crystal over $\Spec (\mathbb F_p)$ of rank $1$ and Newton polygon slope $b$ defined by the pair $(\mathbb Z_p,p^b1_{\mathbb Z_p})$. Let $\eta$ be the generic point $\Spec K\to S$ of $S$. Let $s,l\in\mathbb N^{\ast}$.

In Subsection \ref{S31} we consider commutative affine group schemes $\mathbb H_s$ over $S$ of morphisms between certain evaluations of $\mathcal E_b$ and $\mathcal C$. In Sections \ref{S32} we glue morphisms between different such evaluations in order to introduce good sections above $T$ of the morphisms $\mathbb H_s\to S$ in Subsection \ref{S33}. In Subsection \ref{S34} we complete the proof of Theorem \ref{T1}. The key idea (the plan) can be summarized as follows: under suitable reductions, for $s>>0$ via such good sections above $T$ we can identify $T$ with a closed subscheme of $\mathbb H_s$ and therefore we can conclude that $T$ is an affine scheme.

If $R$ is a reduced perfect ring of characteristic $p$, following \cite{K} we say that an $F^n$-crystal $\mathcal F$ over $\Spec R$ is divisible by $b$ if its evaluation at the endomorphism $\Phi_R^n$ of the thickening $\mathbb W(R)$ is defined by a $\Phi_R^n$-linear endomorphism whose $q$-th iterate for all $q\in\mathbb N^{\ast}$ is congruent to $0$ modulo $p^{bq}$. Thus if $y\in\Spec R$, then the Hodge polygon slopes of $\mathcal F_y$ are all greater or equal to $b$. 

\subsection{Moduli group schemes of morphisms}\label{S31}

For an $A$-algebra $B$ and an $F^n$-crystal $\mathcal F$ over $B$, let $\mathbb E_s(\mathcal F)$ be the evaluation of $\mathcal F$ at the thickening $\mathbb W_s(B)$; it is an object of the category $\mathcal M(W_s(B))$. In particular, we write $\mathbb E_s(\mathcal C_B)=(M_{s,B},\phi_{M_{s,B}})$ and let $\mathbb E_s(\mathcal E_{b,B})=(N_{s,B},\phi_{N_{s,B}})$. Thus we have $M=M_{1,A}$, $\phi_M=\phi_{M_{1,A}}$, and $N_{s,B}=W_s(B)$. Moreover $\phi_{N_{s,B}}: N_{s,B}\to N_{s,B}$ is the $\Phi_B^n$-linear endomorphism which maps $1$ to $p^b$ and $\phi_{M_{s,B}}: M_{s,B}\to M_{s,B}$ is a $\Phi_B^n$-linear endomorphism and we have $M_{s,B}=W_s(B)\otimes_{W_s(A)} M_{s,A}$. The kernel of the epimorphism $W_s(B)\to W_1(B)=B$ is a nilpotent ideal. Based on this and the fact that $M$ is a free $A$-module of rank $r$, we get that each $M_{s,B}$ is a free $W_s(B)$-module of rank $r$. 

We consider the commutative affine group scheme $\mathbb H_s$ over $S$ which represents the following functor: for an $A$-algebra $B$, the abelian group
$$\mathbb H_s(B)=\Hom_{\mathcal M(W_s(B))}(\mathbb E_s(\mathcal E_{b,B}),\mathbb E_s(\mathcal C_{B}))$$
is the group of all $W_s(B)$-linear maps $f:N_{s,B}\to M_{s,B}$ which satisfy the identity $f\circ \phi_{N_{s,B}}=\phi_{M_{s,B}}\circ f$. The $S$-scheme $\mathbb H_s$ is of finite presentation (for $n=1$, see \cite{V1}, Lemma 2.8.4.1; the proof of loc. cit. applies to all $n\in\mathbb N^{\ast}$). 

Let $x\in S$ be a point of codimension $1$. Thus the local ring $D_x:=\mathcal O_{S,x}$ of $S$ at $x$ is a discrete valuation ring. Let $E_x$ be a complete discrete valuation ring which dominates $D_x$ and has a residue field which is algebraically closed. Let $P_x$ be the perfection of $E_x$. We recall that $\mathcal C_{P_x}$ is the pull back of $\mathcal C$ via the natural morphism $\Spec P_x\to S$. As $S=S_{\ge\nu_1}$, the Newton polygon slopes of the two fibers of $\mathcal C_{P_x}$ are greater or equal to $b$. Thus from \cite{K}, Theorem 2.6.1 we get the existence of an $F^n$-crystal $\mathcal D$ over $\Spec P_x$ which is divisible by $b$ and which is equipped with an isogeny
$$\psi_x:\mathcal D\to \mathcal C_{P_x}$$
whose cokernel is annihilated by $p^t$ for some $t\in\mathbb N$. Based on the proof of loc. cit. we can assume that 
$$t=(r-1)b$$ 
depends only on $r$ and $b$. 

\begin{prop}\label{P1}
We assume that the point $x\in S$ of codimension $1$ belongs to $T$. Then there exists a unique $F^n$-subcrystal $\mathcal D_b$ of $\mathcal D$ which is isomorphic to the pull back $\mathcal E_{b,P_x}$ of $\mathcal E_b$. Moreover, $\mathcal D_b$ has a unique direct supplement in $\mathcal D$.
\end{prop}

\noindent
{\bf Proof:} We know that for $y\in\Spec P_x$ all Hodge polygon slopes of $\mathcal D_y$ are at least $b$. If all Hodge polygon slopes of $\mathcal D_y$ are at least $b+1$, then all Newton polygon slopes of $\mathcal D_y$ are at least $b+1$. As under the morphism $\Spec P_x\to S$, the point $y$ maps to either $x\in T$ or $\eta\in T$ and as $\psi_x$ is an isogeny, $(1,b)$ is a break point of the Newton polygon of $\mathcal D_y$. From the last three sentences we get that $(1,b)$ is a point of the Hodge polygon of $\mathcal D_y$. 

Thus for each point $y\in\Spec P_x$, $(1,b)$ is a break point of the Newton polygon of $\mathcal D_y$ and is a point of the Hodge polygon of $\mathcal D_y$. Due to this, from \cite{K}, Theorem 2.4.2 we get that there exists a unique direct sum decomposition
$$\mathcal D=\mathcal D_b\oplus \mathcal D_{>b}$$
into $F^n$-crystals over $\Spec P_x$, where $\mathcal D_b$ is of rank $1$ and each fiber of it at a point $y\in\Spec P_x$ has all Hodge and Newton polygon slopes equal to $b$ and where $\mathcal D_{>b}$ is of rank $r-1$ and each fiber of it at a point $y\in\Spec P_x$ has all Newton polygon slopes greater than $b$ (and has all Hodge polygon slopes greater or equal to $b$). 

As $\mathcal D$ is divisible by $b$, $\mathcal D_b$ and $\mathcal D_{>b}$ are also divisible by $b$. 

As $P_x$ is perfect, for each $l\in\mathbb N^{\ast}$ we have $W(P_x)/(p^l)=W_l(P_x)$ and the module of differentials $\Omega^1_{W_l(P_x)}$ is $0$. Thus, from \cite{BM}, Proposition 1.3.3 we get that an $F^n$-crystal over $\Spec P_x$ is uniquely determined by its evaluation at the thickening $\mathbb W(P_x)$. The evaluation of $\mathcal E_{b,P_x}$ at the thickening $\mathbb W(P_x)$ is canonically identified with $(W(P_x),p^b\Phi^n_{P_x})$ and the evaluation of $\mathcal D_b$ at the thickening $\mathbb W(P_x)$ can be identified with $(W(P_x),p^b\Phi_b)$, where $\Phi_b: W(P_x)\to W(P_x)$ is a $\Phi_{P_x}^n$-linear endomorphism such that $\Phi_b(1)$ generates $W(P_x)$. 

As $P_x$ is the perfection of $E_x$ and as $E_x$ is complete and has an algebraically closed residue field, the rings $W(P_x)$ and $W_l(P_x)$ are strictly henselian and $p$-adically complete. We check that these properties imply that there exists a unit $\upsilon$ of $W(P_x)$ such that we have 
$$\Phi_b(\upsilon)=\Phi_{P_x}^n(\upsilon)\Phi_b(1)=\upsilon.$$ 
If $n=1$, then from \cite{BM}, Proposition 2.4.9 we get that for each $l\in\mathbb N^{\ast}$ there exists a unit $\upsilon_l\in W(P_x)$ such that we have $\Phi_b(\upsilon_l)-\upsilon_l\in p^lW(P_x)$ and the proof of loc. cit. checks that we can assume that $\upsilon_{l+1}-\upsilon_l\in p^lW(P_x)$. Thus for $n=1$ we can take $\upsilon$ to be the $p$-adic limit of the sequence $(\upsilon_l)_{l\ge 1}$. This argument applies entirely for $n>1$. 

The multiplication by $u$ defines an isomorphism 
$$(W(P_x),p^b\Phi^n_{P_x})\to (W(P_x),p^b\Phi_b)$$
which defines an isomorphim $\mathcal E_{b,P_x}\to\mathcal D_b$.\endproof

\medskip
From now we will assume that $x\in T$. We consider a composite morphism
$$j_x[s]:\mathbb E_s(\mathcal E_{b,P_x})\to \mathbb E_s(\mathcal D_b)\to \mathbb E_s(\mathcal D)=\mathbb E_s(\mathcal D_b)\oplus \mathbb E_s(\mathcal D_{>b})$$
in which the first arrow is an isomorphism and the second arrow is the split monomorphism associated to the direct sum decomposition.

Let $$i_x(s):\mathbb E_s(\mathcal E_{b,P_x})\to \mathbb E_s(\mathcal C_{P_x})$$ be the composite of $j_x[s]$ with the morphism $\psi_x[s]:\mathbb E_s(\mathcal D)\to \mathbb E_s(\mathcal C_{P_x})$ which is the evaluation of the isogeny $\psi_x$ at the thickening $\mathbb W_s(P_x)$ (i.e., which is the reduction modulo $p^s$ of $\psi_x$). From now on, we will take $s>t=(r-1)b$. We note that $\psi_x[s]$ is a quasi-isogeny whose cokernel is annihilated by $p^t$ and whose domain is divisible by $b$.

\subsection{Gluing morphisms}\label{S32}

For each point $x\in T$ of codimension $1$ (i.e., whose local ring $D_x$ is a discrete valuation ring), we follow \cite{V1}, Subsection 2.8.3 to show the existence of a finite field extension $K_x$ of $K$ and of an open subset $T_x$ of the normalization of $T$ in $\Spec K_x$ such that $T_x$ has a local ring which is a discrete valuation ring $D^+_x$ that dominates $D_x$ and moreover we have a morphism 
$$i_{T_x}(s):\mathbb E_s(\mathcal E_{b,T_x})\to \mathbb E_s(\mathcal C_{T_x})$$
of the category $\mathcal M(W_s(T_x))$ which is the composite of a split monomorphism with a quasi-isogeny whose cokernel is annihilated by $p^t$ and whose domain is divisible by $b$. 

To check this, with the notations of Subsection \ref{S31} we consider four identifications $E_s(\mathcal C_{D_x})=(W_s(D_x)^r,\phi_{s,x})$, $\mathbb E_s(\mathcal E_{b,D_x})=(W_s(D_x),p^b\Phi_{D_x}^n)$, $\mathbb E_s(\mathcal D_b)=(W_s(P_x),p^b\Phi_{P_x}^n)$, and $\mathbb E_s(\mathcal D_{>b})=(W_s(P_x)^{r-1},p^b\phi_{s,>b,x})$. Now, the $W_s(P_x)$-linear map $\psi_{s,P_x}:W_s(P_x)^r\to W_s(P_x)^r$ defining $\psi_x[s]$ and the $\Phi_{P_x}^n$-linear map $\phi_{s,>b,x}:W_s(P_x)^{r-1}\to W_s(P_x)^{r-1}$ involve a finite number of coordinates of Witt vectors of length $s$ and therefore are defined over $W_s(B_x)$, where $B_x$ is a finitely generated $D_x$-subalgebra of $P_x$. We can choose $B_x$ such that the resulting $W_s(B_x)$-linear map $\psi_{s,B_x}:W_s(B_x)^r\to W_s(B_x)^r$ has a cokernel annihilated by $p^t$. The faithfully flat morphism $\Spec B_x\to \Spec D_x$ has quasi-sections (cf. \cite{G2}, Corollary (17.16.2)) and therefore there exists a finite field extension $K_x$ of $K$ and a discrete valuation ring $D_x^+$ of the normalization $T$ in $K_x$ which dominates $D_x$ and for which we have a $D_x$-homomorphism $B_x\to D_x^+$. The $W_s(D^+_x)$-linear map $\psi_{s,D^+_x}:W_s(D^+_x)^r\to W_s(D^+_x)^r$ which is the natural tensorization of $\psi_{s,B_x}$ induces (via restriction to the first factor $W_s(D^+_x)$ of $W_s(D^+_x)^r$) a morphism $i_{D^+_x}(s):\mathbb E_s(\mathcal E_{b,D^+_x})\to \mathbb E_s(\mathcal C_{D^+_x})$ of the category $\mathcal M(W_s(D^+_x))$ which is the composite of a split monomorphism with a quasi-isogeny whose cokernel is annihilated by $p^t$ and whose domain is divisible by $b$. It is easy to see that there exists an open subset $T_x$ of the normalization of $T$ in $K_x$ which has $D_x^+$ as a local ring and for which there exists a morphism $i_{T_x}(s):\mathbb E_s(\mathcal E_{b,T_x})\to \mathbb E_s(\mathcal C_{T_x})$ of the category $\mathcal M(W_s(T_x))$ that has all the desired properties and that extends the morphism $i_{D^+_x}(s)$ of the category $\mathcal M(W_s(D^+_x))$.
 
By working with $s+l$ instead of $s$, we can assume that there exists $l\in\mathbb N$, $l>>0$ such that $i_{T_x}(s):\mathbb E_s(\mathcal E_{b,T_x})\to \mathbb E_s(\mathcal C_{T_x})$ is the reduction modulo $p^s$ of a morphism 
$$i_{T_x}(s+l):\mathbb E_{s+l}(\mathcal E_{b,T_x})\to \mathbb E_{s+l}(\mathcal C_{T_x})$$ 
of the category $\mathcal M(W_{s+l}(T_x))$.

Let $I_s$ be the set of morphisms $\mathbb E_s(\mathcal E_{b,\bar K})\to \mathbb E_s(\mathcal C_{\bar K})$ which lift to morphisms $\mathbb E_{s+l}(\mathcal E_{b,\bar K})\to \mathbb E_{s+l}(\mathcal C_{\bar K})$ for some $l>>0$. From \cite{V1}, Theorem 5.1.1 (a) (applied for $l>>0$ which depends only on $b$ and $r$) we get that each element of $I_s$ is the evaluation at the thickening $\mathbb W_s(\bar K)$ of a morphism of $F^n$-crystals $\mathcal E_{b,\bar K}\to \mathcal C_{\bar K}$ (strictly speaking loc. cit. is stated for $n=1$ but its proof works for all $n\in\mathbb N^{\ast}$). This implies that $I_s$ is a finite set whose elements are all pull backs of morphisms of $\mathcal M(W_s(L))$, where $L$ is a suitable finite field extension of $K$ contained in $\bar K$. By replacing $S$ with its normalization in $L$, we can assume that $L=K$. As inside $K_x$ we have an identity $D_x^+\cap K=D_x$, inside $W_s(K_x)$ we have an identity $W_s(D_x^+)\cap W_s(K)=W_s(D_x)$. From the last three sentences we get that the pull back $i_{D_x^+}(s)$ of $i_{T_x}(s)$ to a morphism of $\mathcal M(W_s(D_x^+))$ is the pull back of a morphism of $\mathcal M(W_s(D_x))$. Based on this we can assume that there exists an open subscheme $U_x$ of $T$ which contains $x$ and which has the property that there exists a morphism
$$i_{U_x}(s):\mathbb E_s(\mathcal E_{b,U_x})\to \mathbb E_s(\mathcal C_{U_x})$$
of the category $\mathcal M(W_s(U_x))$ such that $i_{T_x}(s)$ is the pull back of it. 

We consider an identification $\mathcal C_{\bar K}=(Q,\phi_Q)$, where $Q=W(\bar K)^r$ and $\phi_Q:Q\to Q$ is a $\Phi_{\bar K}^n$-linear endomorphism. The Newton polygon $\nu_{\eta}$ of $\mathcal C_{\bar K}$ has the Newton polygon slope $b$ with multiplicity $1$ and therefore there exists a unique non-zero direct summand $Q_b$ of $Q$ such that we have $\phi_Q(Q_b)=p^bQ_b$. The rank of the $W(\bar K)$-module $Q_b$ is $1$. Let $z_b\in Q_b$ be such that $Q_b=W(\bar K)z_b$ and $\phi_Q(z_b)=p^bz_b$; it is unique up to multiplication by units of $W(\mathbb F_{p^n})$.

We have a canonical identification $\mathcal E_{b,\bar K}=(W(\bar K),p^b\Phi_K^n)$. The morphism $\mathbb E_s(\mathcal E_{b,\bar K})\to \mathbb E_s(\mathcal C_{\bar K})$ defined by $i_{T_x}(s)$ is an element of $I_s$ and therefore it is the reduction modulo $p^s$ of a morphism $\lambda_x:(W(\bar K),p^b\Phi_K^n)\to (Q,\phi_Q)$ of $F^n$-crystals over $\bar K$. Clearly $\lambda_x(1)\in Q_b$ and thus there exists a unique element $\tau_x\in  W(\mathbb F_{p^n})$ such that we have 
$$\lambda_x(1)=\tau_x z_b.$$
As $i_{T_x}(s)$ is the composite of a split monomorphism with a quasi-isogeny whose cokernel is annihilated by $p^t$ from Fact \ref{F1} applied with $D=W(\bar K)$ we get that $\tau_x$ modulo $p^{t+1}$ is a non-zero element of $W_{t+1}(\mathbb F_{p^n})$. Therefore we can write $\tau_x=p^{t_x}u_x$, where $u_x\in W(\mathbb F_{p^n})$ is a unit and where $t_x\in\{0,\ldots,t\}$. 

From now on, we will take $s>2t$. We consider the morphism
$$\theta_x:=p^{t-t_x}u_x^{-1}i_{U_x}(s):\mathbb E_s(\mathcal E_{b,U_x})\to \mathbb E_s(\mathcal C_{U_x})$$ 
of the category $\mathcal M(W_s(U_x))$; its pull back to a morphism of $\mathcal M(W_s(T_x))$ is the composite of a split monomorphism with a quasi-isogeny whose cokernel is annihilated by $p^{t+t_x}$ and thus also by $p^{2t}$ and whose domain is divisible by $b$. The pull back of $\theta_x$ to a morphism of $\mathcal M(W_s(\bar K))$ is the reduction modulo $p^s$ of the morphism $p^{t-t_x}u_x^{-1}\lambda_x:(W(\bar K),p^b\Phi_K^n)\to (Q,\phi_Q)$ which maps $1$ to $p^tz_b$ and which does not depend on the point $x\in T$ of codimension $1$. 

Let $U$ be the open subscheme of $T$ which is the union of all $U_x$'s. From the previous paragraph we get that the $\theta_x$'s glue together to define a morphism
$$\theta:\mathbb E_s(\mathcal E_{b,U})\to \mathbb E_s(\mathcal C_{U})$$
of the category $\mathcal M(W_s(U))$. 

By replacing $S$ with its normalization in anyone of the finite field extensions $K_x$ of $K$, we can assume that there exists an open dense subscheme $U_0$ of $U$ such that the pull back $\theta_{U_0}:\mathbb E_s(\mathcal E_{b,U_0})\to \mathbb E_s(\mathcal C_{U_0})$ of $\theta$ to a morphism of $\mathcal M(W_s(U_0))$ is the composite of a split monomorphism with a quasi-isogeny whose cokernel is annihilated by $p^{2t}$ and whose domain is divisible by $b$: under such a replacement, we can take $U_0$ to be $T_x$ itself.

\subsection{Good section of $\mathbb H_s$}\label{S33}

We have $\text{codim}_T(T-U)\ge 2$ and the morphism $\theta$ is defined by a section $\theta:U\to \mathbb H_s$ denoted in the same way.

Let $\mathbb I_s$ be the schematic closure $\overline{\theta(U)}$ of $\theta(U)$ in $\mathbb H_s$. As the scheme $\mathbb H_s$ is affine and noetherian and as $U$ is an integral scheme, the scheme $\mathbb I_s$ is also affine, noetherian, and integral. We have a commutative diagram:

 \begin{center}
$\xymatrix{&& \mathbb I_s\ar[d]^{\mbox{affine}}\\
 U\ar@{^{(}->}[r]\ar[urr]^{\mbox{open} \;\;\theta}&T\ar@{^{(}->}[r]&S.}$
\end{center}

We consider the pullback $\mathbb J_s$ of $\mathbb I_s$ to $T$:

 \begin{center}
$\xymatrix{&\mathbb J_s\ar@{^{(}->}[r]^{\mbox{open}}\ar[d]_{\xi}^{\mbox{affine}}& \mathbb I_s\ar[d]^{\mbox{affine}}\\
 U\ar@{^{(}->}_{\mbox{open}}[r]\ar@{^{(}->}[ur]^{\mbox{open}}&T\ar@{^{(}->}[r]_{\mbox{open}}&S.}$
\end{center}

\begin{lemma}\label{L2}
The affine morphism $\xi:\mathbb J_s\to T$ is an isomorphism.
\end{lemma}

\noindent
{\bf Proof:} To prove that $\xi$ is an isomorphism, we can assume that $T=S=\Spec A$ is an affine scheme. As $\xi$ is an affine morphism, $\mathbb J_s=\Spec B$ is also an affine scheme. Since $U$ is open dense in both $T$ and $\mathbb I_s$, $T$ and $\mathbb J_s$ have the same field of fractions $K$. As $\mbox{codim}_T(T-U)\ge 2$ and as $U$ is an open subscheme of both $T$ and $\mathbb J_s$, we have $A_{\mathfrak p}=B_{\mathfrak p}$ for each prime $\mathfrak p\in S=T$ of height $1$. As $A$ is a noetherian normal domain, inside $K$ we have 
$$A\subseteq B\subseteq\displaystyle\bigcap_{\mathfrak q\in\Spec B\mbox{ of height }1}B_{\mathfrak q}\subseteq \displaystyle\bigcap_{\mathfrak p\in \Spec A\mbox{ of height }1}A_{\mathfrak p}=A$$ (cf. \cite{M}, (17.H), Theorem 38 for the equality part; the first inclusion is defined by $\xi$). Therefore $A=B$.\endproof

\medskip
This Lemma \ref{L2} allows us in what follows to identify $T$ itself with an open dense subscheme of $\mathbb I_s$ (i.e., with $\mathbb J_s$).

\subsection{End of the proof}\label{S34}

In this subsection we will show that for $s>>0$, we have $T=\mathbb I_s$. This will complete the proof of Theorem \ref{T1} as $\mathbb I_s$ is an affine scheme. 

We are left to show that the assumption that for $s>>0$ we have $T\neq \mathbb I_s$ leads to a contradiction. This assumption implies that there exists an algebraically closed field $k$ of characteristic $p$ and a morphism $\zeta_0:\Spec (k[[X]])\to \mathbb I_s$ with the properties that under it the generic point of $\Spec (k[[X]])$ maps to $U_0$ and its special point maps to $\mathbb I_s-T$. 

Let $P=k[[X]]^{\perf}$ be the perfection of $k[[X]]$, let $\kappa$ be the perfect field which is the field of fractions of $P$, and let $\zeta:\Spec P\to \mathbb I_s$ be the morphism defined naturally by $\zeta_0$. To the composite of $\zeta$ with the closed embedding $\mathbb I_s\to\mathbb H_s$ corresponds a morphism
$$\omega:\mathbb E_s(\mathcal E_{b,P})\to \mathbb E_s(\mathcal C_P)$$
of the category $\mathcal M(W_s(P))$ whose pull back $\omega_{\kappa}$ to a morphism of $\mathcal M(W_s(\kappa))$ is equal to the pull back $\theta_{\kappa}:\mathbb E_s(\mathcal E_{b,\kappa})\to \mathbb E_s(\mathcal C_{\kappa})$ of $\theta$. 

We have a natural identification $\mathbb E_s(\mathcal E_{b,P})=(W_s(P),p^b\Phi_P^n)$ and we consider an identification $\mathbb E_s(\mathcal C_{P})=(W_s(P)^r,\phi)$. Thus we have a $W_s(P)$-linear map
$$\omega:W_s(P)\to W_s(P)^r$$ 
such that $\omega\circ p^b\Phi_P^n=\phi\circ\omega$.
We consider an isogeny $\mathcal D\to \mathcal C_P$ whose cokernel is annihilated by $p^t$ and with $\mathcal D$ divisible by $b$, again cf. \cite{K}, Theorem 2.6.1 (here $t=(r-1)b$ is as before Proposition \ref{P1}). Thus we also have an isogeny $\iota:\mathcal C_P\to\mathcal D$ whose cokernel is annihilated by $p^t$. We consider its evaluation
$$\iota[s]:\mathbb E_s(\mathcal C_P)\to \mathbb E_s(\mathcal D)$$
at the thickening $\mathbb W_s(P)$. Under an identification $\mathbb E_s(\mathcal D)=(W_s(P)^r,p^b\varphi)$ with $\varphi:W_s(P)^r\to W_s(P)^r$ as a $\Phi_P^n$-linear endomorphism, we get a $W_s(P)$-linear endomorphism $\iota[s]:W_s(P)^r\to W_s(P)^r$ such that we have $\iota [s]\circ\phi=p^b\varphi\circ\iota [s]$. We consider the composite morphism
$$\rho=\iota[s]\circ \omega:\mathbb E_s(\mathcal E_{b,P})\to\mathbb E_s(\mathcal D)$$
identified with a $W_s(P)$-linear map $\rho:W_s(P)\to W_s(P)^r$ such that we have $\rho\circ p^b\Phi_P^n=p^b\varphi\circ\rho$. Let 
$$\gamma=\rho(1)=(\gamma_1,\ldots,\gamma_r)\in W_s(P)^r.$$
From the identity $\rho\circ p^b\Phi_P=p^b\varphi\circ\rho$ we get that the image of $\varphi(\gamma)-\gamma$ in $W_{s-b}(P)^r$ is $0$.
Writing $\gamma=p^u\delta$, where $u\in\mathbb N$ and $\delta\in W_s(P)^r-pW_s(P)^r$, we get that the image of $\varphi(\delta)-\delta$ in $W_{s-b-u}(P)^r$ is $0$. Let $\bar\delta\in P^r-0$ be the image in $P^r=W_1(P)$ of $\delta$ (i.e., the reduction modulo $p$ of $\delta$). 

\begin{lemma}\label{L3}
If $s\ge 3t+1$, then we have $u\le 3t$. Therefore, if moreover we have $s\ge 3t+b+1$, then the image of $\varphi(\delta)-\delta$ in $W_{s-b-3t}(P)^r$ is $0$.
\end{lemma}

\noindent
{\bf Proof:} To check this we can work over $W_s(\kappa)$. As the generic point of $\Spec P$ maps to $U_0$, $\omega_{\kappa}=\theta_{\kappa}:\mathbb E_s(\mathcal E_{b,\kappa})\to \mathbb E_s(\mathcal C_{\kappa})$ is the pull back of the morphism $\theta_{U_0}$. The pull back $\rho_{\kappa}$ of $\rho$ to $\mathcal M(W_s(\kappa))$ is a composite morphism
$$\rho_{\kappa}=\iota[s]_{\kappa}\circ \theta_{\kappa}:\mathbb E_s(\mathcal E_{b,\kappa})\to \mathbb E_s(\mathcal D_{\kappa})$$
and therefore it is the composite of a split monomorphism with a quasi-isogeny whose cokernel is annihilated by $p^{2t}$ (as $\theta_{U_0}$ has this property) and with a quasi-isogeny whose cokernel is annihilated by $p^t$ (as $\iota$ is an isogeny whose cokernel is annihilated by $p^t$). Therefore, $\rho_{\kappa}$ is also the composite of a split monomorphism with a quasi-isogeny whose cokernel is annihilated by $p^{3t}$. This implies that the image of $\gamma$ in $W_{3t+1}(\kappa)$ is non-zero (cf. Fact \ref{F1} applied with $D=W(\kappa)$) and therefore we have $u\le 3t$.\endproof

\begin{lemma}\label{L4}
If $s\ge 3t+b+1$, then the image of $\bar\delta$ in $k^r=W_1(k)^r$ is non-zero. 
\end{lemma}

\noindent
{\bf Proof:} We show that the assumption that the image of $\bar\delta\in P^r-0$ in $k^r=W_1(k)^r$ is $0$ leads to a contradiction. This assumption implies that there exists a largest positive number $c$ of denominator a power of $p$ such that we have 
$$\bar\delta\in X^cP^r\subset P^r=(k[[X]]^{\perf})^r.$$ 
Let $\bar\varphi:P^r\to P^r$ be the $P$-linear endomorphism which is the reduction modulo $p$ of $\varphi$. From Lemma \ref{L3} we get that $\bar\delta=\bar\varphi(\bar\delta)$. Thus
$\bar\delta\in\bar\varphi(X^cP^r)\subseteq X^{p^nc}P^r$
and this implies that $p^nc\le c$ which is a contradiction.
\endproof

\medskip
From the inequality $u\le 3t$ (see Lemma \ref{L3}) and from Lemma \ref{L4} we get that for $s\ge 3t+b+1$ the pull back $\omega_k$ of $\omega$ to a morphism of $\mathcal M(W_s(k))$ is such that its reduction modulo $p^{3t+1}$ is non-zero. For $s>3t+b+1+l$ with $l\in\mathbb N^{\ast}$ large enough but depending only on $b$ and $r$, the reduction of $\omega_k$ modulo $p^{s-l}$ lifts to a morphism $\mathcal E_{0,k}\to \mathcal D_k$ (cf. \cite{V1}, Theorem 5.1.1 (a); again loc. cit. stated for $n=1$ applies to all $n\in\mathbb N^*$) which is non-zero. Thus $\mathcal D_k$ has Newton polygon slope $b$ with multiplicity at least $1$. From this and the existence of the isogeny $\iota$ we get that $\mathcal C_k$ has Newton polygon slope $b$ with multiplicity at least $1$. This implies that the special point of $\Spec (k[[X]])$ under the composite of $\zeta_0:\Spec (k[[X]])\to \mathbb I_s$ with the morphism $\mathbb I_s\to S$ does not map to a point of $S_{\nu_2}=S-T$ and therefore it maps to a point of $T$. Contradiction. This ends the proof of Theorem \ref{T1}.\endproof

\section{Applications of Theorem \ref{T1}}\label{S4}

In Subsection \ref{S41} we prove Corollary \ref{C1}. In Subsection \ref{S42} we follow \cite{V2} to introduce generalized Artin--Schreier systems of equations and Artin--Schreier stratifications. In Subsection \ref{S43} we refine and reobtain Corollary \ref{C2} in the context of these stratifications. Subsection \ref{S44} contains some complements, including Proposition \ref{P3} which prove that `pure in' implies `weakly pure in'. Until the end let $A$ be an arbitrary $\mathbb F_p$-algebra. 

\subsection{Proof of Corollary \ref{C1}}\label{S41}

To prove Corollary \ref{C2}, in this subsection we can assume that $S=\Spec A$ and $d\in\mathbb N$ are as in the paragraph before Subsection \ref{S21}. We can also assume that $\nu(r)=d$ as otherwise $S_{\nu}=\emptyset$ is pure in $S$. Let $l\in\mathbb N$ be such that the Newton polygon $\nu$ has exactly $l+1$ breaking points denoted as $(a_0,b_0)=(0,0), \ldots, (a_l,b_l)=(r,d)$. 

We have obvious identities
$$S_{\nu}=[S_{\ge \nu}\bigcap_{i=0}^l T_{(a_l,b_l)}(\mathcal C)]_{\text{red}}=[S_{\ge \nu}\times_S (T_{(a_0,b_0)}(\mathcal C))_S\times\cdots\times_S T_{(a_l,b_l)}(\mathcal C)]_{\text{red}}.$$
From Theorem \ref{T1} we get that each $T_{(a_l,b_l)}(\mathcal C)$ is an affine scheme. We recall that $S_{\ge\nu}$ is a reduced closed subscheme of $S$. From the last three sentences we get that $S_{\nu}$ is an affine scheme, i.e., is pure in $S$.\endproof  

\subsection{Artin--Schreier stratifications}\label{S42}

Let $x_0$, $x_1$,\ldots, $x_r$ be free variables. For $i,j\in\{1,\ldots,r\}$ let $P_{i,j}(x_0)\in A[x_0]$ be a polynomial which is a linear combination with coefficients in $A$ of the monomials $x_0^q$ with $q\in\mathbb N$ either $0$ or a power of $p$. By a generalized Artin--Schreier system of equations in $r$ variables over $A$ we mean a system of equations of the form
$$x_i=\sum_{j=1}^r P_{i,j}(x_j^p)\;\;\; i\in\{1,\ldots,r\}$$
to which we associate the $A$-algebra 
$$B=A[x_1,\ldots,x_r]/(x_1-\sum_{j=1}^r P_{1,j}(x_j^p),x_2-\sum_{j=1}^r P_{2,j}(x_j^p),\ldots,x_r-\sum_{j=1}^r P_{r,j}(x_j^p)).$$ 
Each equation of the form $x_i=\sum_{j=1}^r P_{i,j}(x_j^p)$ will be called as a generalized Artin--Schreier equation, and its degree $e_i\in\mathbb N$ is defined as follows. We have $e_i=0$ if and only if for all $j\in\{1,\ldots,r\}$ the polynomial $P_{i,j}(x_0)$ is a constant, and if $e_i>0$ then $e_i$ is the largest integer such that there exists a $j\in\{1,\ldots,r\}$ with the property that the degree of $P_{i,j}(x_j^p)$ is $p^{e_i}$. 

Let $e=\text{max}\{e_1,\ldots,e_r\}$; we call it the degree of the generalized Artin--Schreier system of equations in $r$ variables over $A$. Following \cite{V2}, when $e\le 1$ we drop the word `generalized'. 

\begin{prop}\label{P2}
The morphism $\epsilon: \Spec B\to\Spec A$ is \'etale and surjective and its geometric fibers have a number of points equal to a power of $p$. 
\end{prop}

\noindent
{\bf Proof:} If $e_i>1$, then by adding for each $j\in\{1,\ldots,r\}$ such that the degree of $P_{i,j}(x_j^p)$ is $p^{e_i}$ an extra variable $y_{i,j}$ and an equation of the form $y_{i,j}=x_j^p$, the generalized Artin--Schreier equation $x_i=\sum_{j=1}^r P_{i,j}(x_j^p)$ gets replaced by several generalized Artin--Schreier equations of degrees less than $e_i$. By repeating this process of adding extra variables and equations which (up to isomorphisms between $\Spec A$-schemes) does not change the morphism $\epsilon:\Spec B\to\Spec A$, we can assume that $e\le 1$. Thus the proposition follows from \cite{V2}, Theorem 2.4.1 (a) and (b).\endproof

\medskip
The below definition is a natural extrapolation of \cite{V2}, Definition 2.4.2 which applies to \'etale morphisms $\epsilon: \Spec B\to\Spec A$ as in Proposition \ref{P2}.

\begin{df}\label{D1}
Let $\varepsilon:\Spec\mathcal B\to\Spec A$ be an \'etale morphism between affine $\mathbb F_p$-schemes. 

\medskip
{\bf (a)} We assume that $A$ is noetherian. Then by the Artin--Schreier stratification of $\Spec A$ associated to $\varepsilon: \Spec\mathcal B\to\Spec A$ in reduced locally closed subschemes $\mathcal V_1, \ldots, V_q$ we mean the stratification defined inductively by the following property: for each $l\in\{1,\ldots,q\}$ the scheme $V_l$ is the maximal open subscheme of the reduced scheme of $(\Spec A)-(\cup_{q=1}^{l-1} V_q)$ which has the property that the morphism $\epsilon_{V_l}:(\Spec B)\times_{\Spec A} V_l\to V_l$ is an \'etale cover.

\smallskip
{\bf (b)} Let $\mu_1>\mu_2>\cdots >\mu_v$ be the shortest sequence of strictly decreasing natural numbers such that each fiber of the morphism $\epsilon: \Spec B\to\Spec A$ has a number of geometric points equal to $\mu_l$ for some $l\in\{1,\ldots,v\}$. Then by the functorial Artin--Schreier stratification of $\Spec A$ associated to $\varepsilon: \Spec\mathcal B\to\Spec A$ we mean the stratification of $\Spec A$ in reduced locally closed subschemes $U_1, \ldots, U_v$ defined inductively by the following property: for each $l\in\{1,\ldots,v\}$ the scheme $U_l$ is the maximal open subscheme of the reduced scheme of $(\Spec A)-(\cup_{q=1}^{l-1} U_q)$ which has the property that the morphism $\epsilon_{U_l}:(\Spec B)\times_{\Spec A} U_l\to U_l$ is an \'etale cover whose all fibers have a number of geometric points equal to $\mu_l$.
\end{df}

The existence of the stratification $V_1, \ldots, V_q$ of $\Spec A$ is a standard piece of algebraic geometry. The existence of the sequence $\mu_1>\mu_2>\cdots >\mu_v$ follows from the facts that  each \'etale morphism is locally quasi-finite and that $\Spec\mathcal B$ is quasi-compact. The existence of the stratification $U_1, \ldots, U_v$ of $\Spec A$ is implied by  \cite{G2}, Proposition 18.2.8 and Corollary 18.2.9 which show that one can define $U_l$ directly and functorially as follows: each $U_l$ is the set of all points $x\in\Spec A$ such that the fiber of $\varepsilon$ at $x$ has exactly $\mu_l$ geometric points. 

\begin{thm}\label{T2}
Let $\varepsilon:\Spec\mathcal B\to\Spec A$ be an \'etale morphism between affine $\mathbb F_p$-schemes. Then the functorial Artin--Schreier stratification of $\Spec A$ associated to $\varepsilon: \Spec\mathcal B\to\Spec A$ in reduced locally closed subschemes $U_1, \ldots, U_v$ is pure, i.e., for each $l\in\{1,\ldots,v\}$ the stratum $U_l$ is pure in $\Spec A$.
\end{thm}

\noindent
{\bf Proof:} As the \'etale morphism $\varepsilon:\Spec\mathcal B\to\Spec A$ is of finite presentation and due to the functorial part, we can assume that $A$ is a finitely generated $\mathbb F_p$-algebra and thus an excellent ring. We follow \cite{V3}. By replacing $\Spec A$ by its closed subscheme $(\Spec A)-(\cup_{q=1}^{l-1} U_q)$ endowed with the reduced structure, we can assume that $l=1$ and that $A$ is reduced. Thus $U_1$ is an open dense subscheme of $\Spec A$. Based again on \cite{V1}, Lemma 2.9.2 to prove that $U_1$ is an affine scheme, we can replace $A$ by its normalization in its ring of fractions. Thus by passing to connected components of $\Spec A$, we can assume that $A$ is an excellent normal domain. Thus $B=\prod_{l=1}^w B_l$ is a finite product of excellent normal domain which are \'etale $A$-algebras. Let $K_l$ be the field of fractions of $B_l$. Let $L$ be the finite Galois extension of the field of fractions $K$ of $A$ generated by the finite separable extensions $K_l$'s of $K$. By replacing $A$ by its normalization in $L$ (again based on \cite{V1}, Lemma 2.9.2), we can assume that we have $K=K_1=\cdots=K_w$. This implies that each $\Spec (B_l)$ is an open subscheme of $\Spec A$ and thus 
$$U_1=\cap_{l=1}^w \Spec (B_l)=(\Spec (B_1))\times_{\Spec A} (\Spec (B_2))\times_{\Spec A} \cdots\times_{\Spec A} (\Spec (B_w))$$
is the affine scheme $\Spec (B_1\otimes_A\otimes\cdots\otimes_A B_w)$.\endproof

\subsection{A second proof of Corollary \ref{C2}}\label{S43}

We will use Theorem \ref{T2} to obtain a second proof of Corollary \ref{C2} which is simpler and independent of Theorem \ref{T1}. We can assume that $S=\Spec A$ is affine and let $\phi_M:M\to M$ be as in Subsection \ref{S21}.

The identities $S_m=T_{(m,0)}(\mathcal C)$ if $m>0$ and $S_0=T_{(1,0)}(\mathcal C\oplus \mathcal E_0)$ show that $S_m$ is a reduced locally closed subscheme of $S$. Thus by replacing $S$ by $\bar S_m$, we can assume that $S_m$ is an open dense subscheme of $S=\bar S_m$.

We consider the equation
\begin{equation}\label{E1}
\phi_M(z)=z
\end{equation}
in $z\in M$. For $x\in S$ we have $\chi(x)=\dim_{\mathbb F_{p^n}}(\vartheta_x)$, where $\vartheta_x$ is the $\mathbb F_{p^n}$-vector space of solutions of the tensorization of the Equation \ref{E1} over $A$ with an algebraic closure of the residue field $k(x)$ of $S$ at $x$.

From now on we will forget about $\mathcal C$ and just work with the free $A$-module $M$ of rank $r$ and its $\Phi_A^n$-linear endomorphism $\phi_M:M\to M$ and we only assume that we have an open dense subset $S_m$ of $S=\Spec A$ defined by the following property: for $x\in S$, we have $x\in S_m$ if and only if $\dim_{\mathbb F_{p^n}}(\vartheta_x)=m$.

With respect to a fixed $A$-basis $\{v_1,\ldots,v_r\}$ of $M$, by writing $z=\sum_{i=1}^r x_iv_i$ the Equation \ref{E1} defines a generalized Artin--Schreier system of equations in the $r$ variables $x_1,\ldots, x_r$ of the form 
$$x_i=L_i(x_1^{p^n},\ldots,x_r^{p^n})\;\;\; i\in\{1,\ldots,r\},$$
where each $L_i$ is a homogeneous polynomial of total degree at most $1$. Let
$$B=A[x_1,\ldots,x_r]/(x_1-L_1(x_1^{p^n},\ldots,x_r^{p^n}),\ldots,x_r-L_r(x_1^{p^n},\ldots,x_r^{p^n})),$$ 
let $\epsilon:\Spec B\to S$ and let $U_1,\ldots,U_v$ be the functorial Artin--Schreier stratification of $S$ associated to $\epsilon: \Spec B\to S$. Let $p^{\mu_1}>p^{\mu_2}>\cdots >p^{\mu_v}$ be the shortest sequence of strictly decreasing of powers of $p$ by natural numbers such that for each each $l\in\{1,\ldots,v\}$ every geometric fiber of the morphism $\epsilon_{U_l}: \Spec B\times_S U_l\to U_l$ has a number of geometric points equal to $p^{\mu_l}$, cf. Proposition \ref{P2} and Definition \ref{D1} (b).

The fact that the morphism $\epsilon:\Spec B\to S$ is \'etale (cf. Proposition \ref{P2}) is equivalent to \cite{Z}, Proposition 3. We consider the lower semi-continuous function (cf. \cite{G2}, Proposition 18.2.8)
$$\mu:S\to\mathbb N$$ 
defined by the rule: $\mu(x)=p^{n\dim_{\mathbb F_{p^n}}(\vartheta_x)}$ is the number of geometric points of $\epsilon:\Spec B\to S$ above $x$ (i.e., is the number of elements of $\vartheta_x$). We get that $\mu_l$ is divisible by $n$ for all $l\in\{1,\ldots,v\}$ and (as $S_m$ is dense in $S$) we have $\mu_1=mn$. Moreover, for $x\in S$ and $q\in\mathbb N$ we have $\mu(x)=p^{nq}$ if and only if $x\in S_q$. We conclude that $S_m=U_1$ and therefore (cf. Theorem \ref{T2}) $S_m$ is an affine scheme.\endproof

\subsection{Complements}\label{S44}

For the sake of completeness, we include a proof of the following well-known result (to be compare with \cite{V1}, Remark 6.3 (a)).

\begin{prop}\label{P3}
Let $Z$ be a reduced locally closed subscheme of a locally noetherian scheme $Y$. If $Z$ is pure in $Y$, then $Z$ is weakly pure in $Y$.
\end{prop}

\noindent
{\bf Proof:} We can assume that $Z\subsetneq \bar Z=Y$. By localizing $Y$ at the generic point of an irreducible component of $\bar Z-Z$, we can assume that $Y=\bar Z=\Spec C$ is a local affine scheme of dimension at least $1$ and $Z$ is the complement in $Y$ of the closed point of $Y$ and we have to prove that $C$ has dimension $1$. By passing to a connected component of the normalization of the reduced completion $\hat C_{\text{red}}$ of $C$ in the ring of fractions of $\hat{C}_{\text{red}}$, we can assume that $C$ is in fact an integral normal local ring which is not a field. 

We show that the assumption that $\dim(C)\ge 2$ leads to a contradiction. As the open dense subscheme $Z$ of $Y$ is pure in in $Y$, $Z$ is the spectrum of a $C$-subalgebra of the field of fractions of $C$ which contains $C$ and which is contained in the intersection of all the localizations of $C$ at points of $Y$ of codimension $1$ in $Y$ (as these points belong to $Z$). As $\dim(C)\ge 2$, from \cite{M}, (17H), Theorem 38 we get that this intersection is $C$ and thus we have $Z=\Spec C=Y$. Contradiction. Thus $\dim(C)=1$.\endproof

\begin{remark}\label{R1}
Suppose $A$ is a local noetherian $\mathbb F_p$-algebra of dimension at least $2$. Let $\mathfrak m$ be the maximal ideal of $A$. Suppose $M=A^r$ is equipped with a $\Phi_A^n$-linear endomorphism $\phi_M:M\to M$ such that for each non-closed point $x$ of $S=\Spec A$, with the notations of Subsection \ref{S43} we have $\dim_{\mathbb F_{p^n}}(\vartheta_x)=m$. Then $S_m=U_1$ being pure in $S$, it is also weekly pure in $S$ (cf. Proposition \ref{P3}) and thus $S-S_m$ cannot be $\mathfrak m$ as $\text{codim}_S({\mathfrak m})\ge 2$. Therefore we have $S_m=S$ and in this way we reobtain \cite{Z}, Proposition 5. One can view Theorem \ref{T2} as a generalization and a refinement of \cite{Z}, Proposition 5.
\end{remark}

\begin{remark}\label{R2}
For $q\in\mathbb N^{\ast}$ we define recursively an $A$-linear map 
$$\phi_M^{(q)}:A\otimes_{F_A^{nq},A} M\to M$$ as follows: let $\phi_M^{(1)}:A\otimes_{F_A^n,A} M\to M$ be the $A$-linear map defined by $\phi_M$ and we have the recursive formula 
$\phi_M^{(q)}=\phi_M^{(1)}\circ (1_A\otimes_{F_A^n, A} \phi_M^{(q-1)}).$
Deligne proved in \cite{D} the case $n=1$ of Theorem \ref{T2} using ranks of images of $\phi_M^{(q)}$ for $q>>0$ at points $x\in S=\Spec A$ and properties of Grassmannians. 
\end{remark}

\bigskip\noindent
{\bf Acknowledgement.}
The first author would like to thank the second author for the continuous support during his Ph.D. studies and his family for spiritual support throughout his life. The second author would like to thank Binghamton University for good working conditions. Both authors would like to thank the referee for many valuable comments.

\hbox{}
\hbox{Jinghao Li,\;\;\;E-mail: jinghao.lee@gmail.com}
\hbox{Address: Sequoia Capital Global Equities
2800 Sand Hill Road, Suite 101,} 
\hbox{Menlo Park, CA 94025, U.S.A.}

\bigskip

\hbox{Adrian Vasiu,\;\;\;E-mail: adrian@math.binghamton.edu}
\hbox{Address: Department of Mathematical Sciences, Binghamton University,}
\hbox{P. O. Box 6000, Binghamton, New York 13902-6000, U.S.A.}
\end{document}